\documentclass[12pt]{article}
\usepackage{amsthm,amsfonts,amssymb,amscd}

\begin{document}


\begin{center}
\large \bf Birationally rigid hypersurfaces \\ with quadratic
singularities of low rank
\end{center}\vspace{0.5cm}

\centerline{A.V.Pukhlikov}\vspace{0.5cm}

\parshape=1
3cm 10cm \noindent {\small \quad\quad\quad \quad\quad\quad\quad
\quad\quad\quad {\bf }\newline It is shown that hypersurfaces of
degree $M$ in ${\mathbb P}^M$, $M\geqslant 5$, with at most
quadratic singularities of rank at least 3, satisfying certain
conditions of general position, are birationally superrigid Fano
varieties and the complement to the set of such hypersurfaces is
for $M\geqslant 8$ of codimension at least ${M-1 \choose 2} + 1$
with respect to the natural parameter space.

Bibliography: 18 items.}

AMS classification: 14E05, 14E07

Key words: Fano variety, birational rigidity, quadratic
singularity.\vspace{0.3cm}

\section*{Introduction}

{\bf 0.1. Statement of the main result.} Let ${\cal P}={\cal
P}_{M,M+1}$ be the linear space of homogeneous polynomials of
degree $M$ in $M+1$ variables, where $M\geqslant 5$, which we
identify with the space of sections $H^0({\mathbb P}^M,{\cal
O}_{{\mathbb P}^M}(M))$. The projectivization of the space ${\cal
P}$ parameterizes hypersurfaces of degree $M$ in ${\mathbb P}^M$.
If for $f\in{\cal P}\backslash\{0\}$ the hypersurface
$F(f)=\{f=0\}$ is irreducible, reduced, factorial and has terminal
singularities, then $F(f)\subset{\mathbb P}^M$ is a primitive Fano
variety:
$$
\mathop{\rm Pic}F(f)={\mathbb Z}H,\quad K_{F(f)}=-H,
$$
where $H$ is the class of a hyperplane section. Let ${\cal F}_{\rm
srigid}\subset{\cal P}$ be the subset, consisting of polynomials
$f\in{\cal P}\backslash\{0\}$, such that the hypersurface $F(f)$
satisfies all properties listed above and, besides, is a
birationally superrigid variety, that is, for every linear system
$\Sigma\subset|nH|$ with no fixed components, where $n\geqslant
1$, and a general divisor $D\in\Sigma$ the pair
$(F(f),\frac{1}{n}D)$ is canonical. Due to the importance of the
property of being birationally (super)rigid for the problems of
higher-dimensional birational geometry (a birationally superrigid
variety $F(f)$ can not be fibred into rationally connected
varieties over a positive-dimensional base, is not birational to
the total space of any Mori fibre space over a
positive-dimensional base, any birational map from this variety
onto a Fano variety of the same dimension with terminal ${\mathbb
Q}$-factorial singularities is an isomorphism), we get the natural
problem of describing, as precisely as possible, the subset ${\cal
F}_{\rm srigid}$ and, in particular, of estimating the codimension
of the complement $\mathop{\rm codim}(({\cal P}\backslash{\cal
F}_{\rm srigid})\subset{\cal P})$.

In \cite{EcklPukhlikov2014} it was shown that ${\cal F}_{\rm
srigid}$ contains the equations of all hypersurfaces with at most
quadratic singularities of rank $\geqslant 5$, satisfying the
point-wise regularity conditions, which implies the estimate
$$
\mathop{\rm codim}(({\cal P}\backslash{\cal F}_{\rm
srigid})\subset{\cal P})\geqslant {{M-3}\choose{2}}+1.
$$

The aim of this paper is to improve this result: we will show that
quadratic singular points of rank 4 and 3, satisfying certain
additional conditions of general position, can also be allowed,
whereas the point-wise regularity conditions can be somewhat
relaxed, which implies a stronger estimate for the codimension of
non-rigid hypersurfaces than the one shown in
\cite{EcklPukhlikov2014}. Let us give precise statements.

Let $f\in{\cal P}\backslash\{0\}$ and $F=F(f)$ be the
corresponding hypersurface. For a point $o\in F$ take an arbitrary
system of affine coordinates $z_1,\dots,z_M$ on an affine chart
${\mathbb A}^M\subset{\mathbb P}^M$, containing that point, where
$o=(0,\dots,0)$, and write down the corresponding affine
polynomial (which for convenience we denote by the same symbol
$f$) as the sum
$$
f(z_1,\dots,z_M)=q_1(z_*)+q_2(z_*)+\dots+q_M(z_*),
$$
where $q_i$ are homogeneous of degree $i$. Obviously,
$o\in\mathop{\rm Sing}F$ if and only if $q_1\equiv 0$. We say that
that a singular point $o\in F$ is {\it quadratic of rank}
$a\geqslant 1$, if $\mathop{\rm rk}q_2=a$.

Assume that the point $o\in F$ is a quadratic singularity of rank
3; in that case we may assume that $q_2=z^2_1+z^2_2+z^2_3$. We say
that the point $o$ satisfies the condition (G), if the singular
locus of the cubic hypersurface
$$
\{q_3(0,0,0,z_4,\dots,z_M)=0\}
$$
in the projective space ${\mathbb P}^{M-4}$ with the homogeneous
coordinates $(z_4:\dots:z_M)$ is zero-dimensional (in particular,
$q_3(0,0,0,z_4,\dots,z_M)\not\equiv 0$) or empty, and moreover,
the homogeneous polynomial
$$
4q_4(0,0,0,z_4,\dots,z_M)-\sum^3_{i=1}\frac{\partial q_3}{\partial
z_i}(0,0,0,z_4,\dots,z_M)^2
$$
of degree 4 does not vanish at any singular point of that
hypersurface.

{\bf Theorem 0.1.} {\it Assume that every point of the
hypersurface $F$ is either non-singular or a quadratic singularity
of rank $\geqslant 3$, and every quadratic singularity of rank 3
satisfies the condition (G) and the inequality
$$
\mathop{\rm codim}(\mathop{\rm Sing}F\subset F)\geqslant 4
$$
holds. Then $F$ is an irreducible reduced factorial variety with
terminal singularities.}

(Of course, the irreducibility, reducedness and factoriality
follow from the inequality for the codimension of the subset
$\mathop{\rm Sing}F$.)

Still using the coordinate notations, introduced above, note that
if $o\in F$ is a non-singular point, then the tangent hyperplane
$T_oF$ is given by the equation $q_1=0$. We say that a
non-singular point $o\in F$ is regular, is for $M\geqslant 6$

(R1) the polynomials
\begin{equation}\label{05.12.23.1}
q_6|_{T_oF},\,\dots,\, q_M|_{T_oF}
\end{equation}
form a regular sequence, that is, the system of equations
$$
q_6|_{T_oF}=\dots=q_M|_{T_oF}=0
$$
has in the hyperplane $T_oF$ a four-dimensional set of common
zeros.

We say that a quadratic singularity $o\in F$ of rank $a\geqslant
7$ (for $M\geqslant 7$) is regular, if

(R2) the polynomials
\begin{equation}\label{05.12.23.2}
q_2,\, q_7,\, \dots,\, q_M
\end{equation}
form a regular sequence, that is, the set of their common zeros is
of dimension 5.

Finally, if $o\in F$ is a quadratic singularity of rank
$\in\{3,4,5,6\}$, we say that it is regular, if

(R3) the polynomials
\begin{equation}\label{05.12.23.3}
q_2,\, q_3,\,\dots,\, q_M
\end{equation}
form a regular sequence, that, the set of their common zeros  is
one-dimensional.

Now let us define the subset ${\cal F}\subset{\cal P}$ if the
following way: $f\in{\cal F}$ if and only if the hypersurface
$F=F(f)$ satisfies the assumptions of Theorem 0.1, and for
$M\geqslant 6$ the inequality
$$
\mathop{\rm codim}(\mathop{\rm Sing}F\subset F)\geqslant 5
$$
holds, and for $M=5$ the hypersurface $F$ does not contain
two-dimensional planes and no section of this hypersurface by a
three-dimensional subspace in ${\mathbb P}^5$ contains a line
consisting of singular point of that section, and (for every
$M\geqslant 5$) every point of $F$ is regular in the sense of the
corresponding condition (R1), (R2) or (R3).

{\bf Theorem 0.2.} {\it For $f\in{\cal F}$ the hypersurface
$F=F(f)$ is a birationally superrigid Fano variety. In particular,
$F$ can not be fibred by a rational map into rationally connected
varieties over a positive-dimensional base, $F$ is not birational
to the total space of any Mori fibre space over a
positive-dimensional base and every birational map $\chi\colon
F\dashrightarrow F'$, where $F'$ is a Fano variety with terminal
${\mathbb Q}$-factorial singularities and the Picard number 1, is
a biregular isomorphism.}

Let us define a function $\gamma\colon \{ M\in {\mathbb Z}\,|\,
M\geqslant 5\}\to {\mathbb Z}$ in the following way:\vspace{0.5cm}

\begin{center}
\begin{tabular}{|c|c|c|c|c|}
\hline
$M$  & 5 & 6 & 7 & $\geqslant 8$ \\
\hline $\gamma(M)$ & 6 & 9 & 15 & ${M-1 \choose 2}+1$ \\
\hline
\end{tabular}
\end{center}
\vspace{0.5cm}

{\bf Theorem 0.3.} {\it The codimension of the complement ${\cal
P}\backslash{\cal F}$ in the space ${\cal P}$ is at least}
$\gamma(M)$.\vspace{0.3cm}


{\bf 0.2. The structure of the paper.} In \S 1 we prove Theorem
0.1. Besides, we describe some local conditions at a quadratic
point $o\in F$ of rank 4 that ensure the inequality $\mathop{\rm
codim}_o(\mathop{\rm Sing} F\subset F)\geqslant 4$ (the symbol
$\mathop{\rm codim}_o$ stands for the codimension in a
neighborhood of the point $o$). We will show that if the condition
(G) is satisfied, then there are finitely many quadratic
singularities of rank 3 and, blowing up this finite set, we get a
variety with at most quadratic points of rank $\geqslant 4$, where
the singular locus of that variety is of codimension $\geqslant
4$.

In \S 2 we prove Theorem 0.2. In order to exclude maximal
singularities over non-singular points, we use the
$8n^2$-inequality \cite{Pukh13a,Ch05c}, in the case of quadratic
singularities of rank $\geqslant 7$ we use the generalized
$4n^2$-inequality \cite{Pukh2017b}, in the case of quadratic
singularities of rank $a\in\{3,4,5,6\}$ we use the traditional
$4n^2$-inequality for quadratic points shown in
\cite{EcklPukhlikov2014} with a correction made in
\cite{Pukh2017a}. In order to obtain a contradiction, we use the
technique of hypertangent divisors based on the regularity
conditions (R1), (R2) and (R3), respectively.

In \S 3 we show Theorem 0.3. Here we use the technique developed
in \cite{Pukh98b,Pukh01,Pukh2017a}.\vspace{0.3cm}


{\bf 0.3. General comments.}  The present paper belongs to the
series of papers on effective birational rigidity, started by
\cite{EcklPukhlikov2014}: for a given family of Fano varieties we
not only prove the birational rigidity of a variety of general
position but also give an explicit estimate for the codimension of
the complement to the set of birationally rigid varieties in that
family. One of the applications of such results is the possibility
to construct fibrations into Fano varieties over a
positive-dimensional base, each fibre of which is a birationally
rigid variety, see the survey \cite{Pukh2022b}.

The approach to studying geometric objects that vary in a certain
family, when at first objects of general position are considered
(in some precise sense, for instance, the objects with no
singularities), and after that more and more special subfamilies
are studied, is sometimes called the ``Poincar\'{e} strategy''. If
we have a family of Fano varieties, defined by means of a
particular geometric construction, then it is natural to study the
problem of their birational rigidity in exactly this way, allowing
more and more complicated singularities. For (irreducible reduced)
factorial hypersurfaces of degree $M$ in ${\mathbb P}^M$ the final
goal is to determine precisely the boundaries of the set of
(super)rigid varieties, in the sense that how ``bad'' the
singularities can be when the variety is still birationally rigid;
for instance, in \cite{Pukh02a} it was shown that the multiplicity
of the unique singular point of general position can be as high as
$M-2$ (and the conditions of general position were relaxed and the
proof simplified in \cite{Pukh2019a}), whereas a hypersurface of
degree $M$ with a point of multiplicity $M-1$ is obviously
rational and can not be birationally rigid.

As an example of such an approach to studying birational geometry
of three-dimensional fibrations into del Pezzo surfaces, when
varieties with more and more complicated singularities are
considered, we refer to the papers
\cite{Krylov2018b,AbbanKrylov22,KrylovOkadaetal22}.

The author is grateful to the members of Divisions of Algebraic
Geometry and Algebra at Steklov Institute of Mathematics for the
interest to his work, and to the colleagues in Algebraic Geometry
research group at the University of Liverpool for general support.


\section{Factorial terminal singularities}

In this section we prove Theorem 0.1. In subsection 1.1 we
consider simple general facts about quadratic singularities, in
Subsection 1.2 we look into what happens with singularities when
we blow up a quadratic point of the given rank. In Subsection 1.3,
on that basis we show Theorem 0.1.\vspace{0.3cm}

{\bf 1.1. Quadratic singularities.} Recall (see
\cite{EcklPukhlikov2014,Pukh15a}), that a point $o\in{\cal X}$ of
an irreducible algebraic variety ${\cal X}$ is a quadratic
singularity of rank $r\geqslant 1$, if in some neighborhood of
that point ${\cal X}$ is realized as a subvariety of a
non-singular $N=(\mathop{\rm dim}{\cal X}+1$)-dimensional variety
${\cal Y}\ni o$, and for some system $(u_1,\dots,u_N)$ of local
parameters on ${\cal Y}$ at the point $o$ the subvariety ${\cal
X}$ is a hypersurface given by an equation
$$
g\in{\cal O}_{o,{\cal Y}}\subset{\mathbb C}[[u_1,\dots,u_N]],
$$
which is represented by the formal power series
$$
g=g_2+g_3+\dots,
$$
where $g_j(u_1,\dots,u_N)$ is a homogeneous polynomial of degree
$j$ and $\mathop{\rm rk}g_2=r$.

The following obvious fact is true.

{\bf Proposition 1.1.} {\it Assume that $o\in{\cal X}$ is a
quadratic singularity of rank $r$. Then in a neighborhood of the
point $o$ every point $p\in{\cal X}$ is either non-singular, or a
quadratic singularity of rank} $\geqslant r$.

Let ${\cal Y}^+\to{\cal Y}$ be the blow up of the point $o$ with
the exceptional divisor $E_{\cal Y}\cong{\mathbb P}^{N-1}$ and
${\cal X}^+\subset{\cal Y}^+$ the strict transform of the
hypersurface ${\cal X}$ on ${\cal Y}^+$, so that ${\cal
X}^+\to{\cal X}$ is the blow up of the point $o$ on ${\cal X}$
with the exceptional divisor $E_{\cal Y}|_{{\cal X}^+}=E_{\cal
X}$. Therefore, $E_{\cal X}$ is a quadric of rank $r$ in the
projective space $E_{\cal Y}$. The following fact is also obvious.

{\bf Proposition 1.2.} {\it The equality $\mathop{\rm
codim}(\mathop{\rm Sing}E_{\cal X}\subset{\cal X})=r$ and the
inequality
$$
\mathop{\rm codim}\nolimits_o(\mathop{\rm Sing} {\cal
X}\subset{\cal X}) \geqslant r-1
$$
hold.}

From here, in view of Grothendieck's theorem on parafactoriality
\cite{CL,Call1994}, it follows immediately that a variety with at
most quadratic points of rank $\geqslant 5$ is factorial (which
was used in \cite{EcklPukhlikov2014}).

Quadratic singularities have the following useful property of
being stable with respect to blow ups, see
\cite{EcklPukhlikov2014} or \cite[Subsection 1.7]{Pukh2022a}.

{\bf Proposition 1.3.} {\it Assume that ${\cal X}$ has at most
quadratic singularities of rank $\geqslant a\geqslant 3$ and
$B\subset{\cal X}$ is an irreducible subvariety of codimension
$\geqslant 2$. Then there is an open subset $U\subset{\cal X}$,
such that $U\cap B\neq\emptyset$, the subvariety $U\subset B$ is
non-singular and its blow up
$$
\sigma_B\colon U_B\to U
$$
gives a quasi-projective variety $U_B$ with at most quadratic
singularities of rank} $\geqslant a$.

(For a proof, based on simple coordinate computations, see
\cite[Subsection 1.7]{Pukh2022a}.)

The arguments given above imply that if ${\cal X}$ is a variety
with at most quadratic singularities of rank $\geqslant 5$, then
${\cal X}$ is a factorial variety with terminal singularities. In
the notations of Proposition 1.3, where $a=5$, for $\mathop{\rm
codim}B\in\{2,3\}$ we may assume that $U$ is non-singular (since
$B\not\subset\mathop{\rm Sing}{\cal X}$), so that the discrepancy
of the exceptional divisor $E_B$ of the blow up $\sigma_B$ is
$\mathop{\rm codim}B-1$, and for $\mathop{\rm codim}B\geqslant 4$
the inequality $a(E_B,U)\geqslant\mathop{\rm codim}B-2$ holds
(where the equality holds precisely when $B\subset\mathop{\rm
Sing}{\cal X}$).

However, in the present paper we allow quadratic singularities of
rank 4 and 3.\vspace{0.3cm}


{\bf 1.2. The blow up of a quadratic point.} Let $o\in{\cal X}$ be
a quadratic singularity of rank $a\geqslant 3$. In the notations
of Subsection 1.1 the local parameters $u_1,\dots,u_N$ can be
identified with homogeneous coordinates on $E_{\cal Y}$, so that
$E_{\cal X}$ is a quadric with the equation $g_2=0$. We may assume
that $g_2=u^2_1+\dots+u^2_a$, so that $\mathop{\rm Sing}E_{\cal
X}=\{u_1=\dots=u_a=0\}$. Set $Q=\{g_3|_{{\rm Sing}E_{\cal
X}}=0\}$, that is, $Q$ is the hypersurface
$g_3(0,\dots,0,u_{a+1},\dots,u_N)=0$ in the projective space
${\mathbb P}^{N-a-1}$ with homogeneous coordinates
$(u_{a+1}:\cdots:u_N)$. (Of course, the linear subspace
$\mathop{\rm Sing}E_{\cal X}$ and the cubic hypersurface
$Q\subset\mathop{\rm Sing}E_{\cal X}$ do not depend on the system
of local parameters at the point $o$.) Set also
$$
h(u_{a+1},\dots,u_N)=4g_4(0,\dots,0,u_{a+1},\dots,u_N)-
\sum^a_{i=1}\frac{\partial g_3}{\partial
u_i}(0,\dots,0,u_{a+1},\dots,u_N)^2.
$$
This is a homogeneous polynomial of degree 4 on the subspace
$\mathop{\rm Sing}E_{\cal X}\subset E_{\cal Y}$.

{\bf Proposition  1.4.} (i) {\it The following equality holds:}
$$
Q=E_{\cal Y}\cap\mathop{\rm Sing}{\cal X}^+.
$$
(ii) {\it A point $p\in Q\backslash\mathop{\rm Sing}Q$ is a
quadratic singularity of rank} $a+2$.

\noindent (iii) {\it A point $p\in\mathop{\rm Sing}Q$ is a
quadratic singularity of rank $a+1$, if $h(p)\neq 0$, and of rank
$a$, otherwise.}

(Note that if $g_3|_{{\rm Sing}E_{\cal X}}\equiv 0$, that is,
$Q=\mathop{\rm Sing}E_{\cal X}$, then $Q=\mathop{\rm Sing}Q$: the
singularities of $Q$ are common zeros of partial derivatives of
the cubic polynomial $g_3(0,\dots,0,u_{a+1},\dots,u_N$).)

{\bf Proof} is obtained by means of obvious local computations.
Let $p\in\mathop{\rm Sing}E_{\cal X}$ be an arbitrary point. We
may assume that
$$
p=(\underbrace{0,\dots,0}\limits_a,1,\underbrace{0,\dots,0}\limits_{N-a-1}),
$$
so that a system of local parameters $w_1,\dots,w_N$ on ${\cal
Y}^+$ at the point $p$ is related to the original system by the
equalities
$$
u_{a+1}=w_{a+1},\,\,u_i=w_iw_{a+1}
$$
for $i\neq a+1$. Now the claim (i) is obvious. If $p\in Q$, then
the first (quadratic) component in the presentation of the local
equation of the variety ${\cal X}^+$ at the point $p$ has the form
$$
\sum^a_{i=1}w^2_i+w_{a+1}\left(\sum_{i\neq
a+1}\alpha_iw_i\right)+w^2_{a+1}g_{4}(0,\dots,0,1,0,\dots,0),
$$
where
$$
\alpha_i=\frac{\partial g_3}{\partial u_i}(0,\dots,0,1,0,\dots,0)
$$
(the only 1 occupies the $(a+1)$-th place). Computing the rank of
this quadratic form, we get the claims (ii) and (iii). Q.E.D. for
the proposition.

If $a=4$, then, since
$$
E_{\cal Y}\cap\mathop{\rm Sing}{\cal X}^+\subset\mathop{\rm
Sing}E_{\cal X},
$$
we obtain the following claim.

{\bf Proposition 1.5.} {\it If $o\in{\cal X}$ is a quadratic
singularity of rank 4 and the inequality $\mathop{\rm
codim}(\mathop{\rm Sing}{\cal X}\subset{\cal X})\geqslant 4$
holds, then the similar inequality is true for ${\cal X}^+$,
either.}

Besides, Proposition 1.4 implies the following fact.

{\bf Proposition 1.6.} {\it If $o\in{\cal X}$ is a quadratic
singularity of rank 4 and either $g_3|_{\mathop{\rm Sing} E_{\cal
X}}\not\equiv 0$, or $g_3|_{\mathop{\rm Sing} E_{\cal X}}\equiv
0$, but $h|_{\mathop{\rm Sing}E_{\cal X}}\not\equiv 0$, then the
inequality}
$$
\mathop{\rm codim}\nolimits_o(\mathop{\rm Sing}{\cal
X}\subset{\cal X})\geqslant 4
$$
{\it holds.}

{\bf Proof.} Indeed, if $Q\neq\mathop{\rm Sing}E_{\cal X}$, then
$$
\mathop{\rm codim}((E_{\cal Y}\cap\mathop{\rm Sing}{\cal
X}^+)\subset E_{\cal X})=4,
$$
which implies the required inequality. If $Q=\mathop{\rm
Sing}E_{\cal X}$, then at the point of general position $p\in Q$
the rank of the quadratic singularity of the variety ${\cal X}^+$
is equal to 5. Therefore, in a neighborhood of the point $o$ the
rank of every singularity of the variety ${\cal X}$ is at least 4,
and the closed set of the points of rank 4 is of codimension
$\geqslant 4$. Now we apply Proposition 1.2 and complete the
proof. Q.E.D.\vspace{0.3cm}


{\bf 1.3. Singularities of rank 3.} Let us prove Theorem 0.1. Let
$o\in F$ be a quadratic point of rank 3, ${\mathbb P}^+\to{\mathbb
P}^M$ the blow up of that point with the exceptional divisor
$E_{\mathbb P}\cong{\mathbb P}^{M-1}$ and $F^+\subset{\mathbb
P}^+$ the strict transform of $F$, so that $E_F=F^+\cap E_{\mathbb
P}$ is a quadric hypersurface of rank 3. In the notations of
Subsection 0.1 set
$$
Q=\{q_3(0,0,0,z_4,\dots,z_M)=0\}=\{q_3|_{\mathop{\rm Sing}
E_F}=0\}
$$
and
$$
h(z_4,\dots,z_M)=4q_4(0,0,0,z_4,\dots,z_M)-\sum^3_{i=1}\frac{\partial
q_3}{\partial z_i}(0,0,0,z_4,\dots,z_M)^2,
$$
these notations agree with Subsection 1.2. No Proposition 1.4
implies that an arbitrary point $p\in E_F$ is either non-singular
on $F^+$ (this happens precisely when $p\not\in Q$), or a
quadratic singularity of rank 5 (when $p\in Q\backslash\mathop{\rm
Sing}Q$), or a quadratic singularity of rank 4 (when
$p\in\mathop{\rm Sing}Q$). In particular, in a neighborhood of the
point $o$ there are no other quadratic singularities of rank 3. We
have shown the following fact.

{\bf Proposition 1.7.} {\it In the assumptions of Theorem 0.1 the
set of quadratic singularities of the hypersurface $F$ of rank 3
is either empty or finite. In the second case, blowing up this
finite set of points, we get a variety that has, at worst,
quadratic singularities of rank $\geqslant 4$, and its singular
locus is of codimension 4.}

Now let $E$ be a prime exceptional divisor over $F$ and $B$ its
centre on $F$, and $o\in B$ a point of general position on $B$. If
$o$ is not a quadratic singularity of rank 3, then Propositions
1.3 and 1.5 imply that $a(E,F)\geqslant 1$. Assume that $B=o$ is a
quadratic singularity of rank 3. Then either $E=E_F$ and
$a(E,F)=M-3\geqslant 2$, or the centre of $E$ on $F^+$ is an
irreducible subvariety $B_+\subset E_F$ of codimension $\geqslant
2$, so that by what was said above about the singularities of the
variety $F^+$ we get
$$
a(E,F)\geqslant a(E_F,F)+1\geqslant M-2.
$$
This proves that the singularities of the variety $F$ are
terminal. Q.E.D. for the proof of Theorem 0.1.


\section{Birational superrigidity}

In this section we prove Theorem 0.2.\vspace{0.3cm}

{\bf 2.1. Types of maximal singularities.} Set $F=F(f)$, where
$f\in{\cal F}$, to be a hypersurface, satisfying the assumptions
of Theorem 0.2. Assume that for some linear system
$\Sigma\subset|nH|$ with no fixed components, where $H=-K_F$ is
the class of a hyperplane section and $n\geqslant 1$, the pair
$(F,\frac{1}{n}D)$, where $D\in\Sigma$ is a general divisor, is
not canonical, that is, for some exceptional divisor $E$ over $F$
(depending on $\Sigma$ only, not on the the divisor $D$) the
Noether-Fano inequality holds:
$$
\mathop{\rm ord}\nolimits_ED>n\cdot a(E),
$$
where $a(E)=a(E,F)$ is the discrepancy of $E$ with respect to $F$.
Let us show that this assumption leads to a contradiction. This
will prove Theorem 0.2.

The exceptional divisor $E$ over $F$ is traditionally called a
{\it maximal singularity} of the mobile linear system $\Sigma$.
Let $B\subset F$ be the centre of the singularity $E$. This is an
irreducible subvariety of codimension $\mathop{\rm codim}(B\subset
F)\geqslant 2$. Depending on the value of this codimension and the
type of a general point $o\in B$, we get the following options:

(1) $\mathop{\rm codim}B\in\{2,3\}$,

(2) $\mathop{\rm codim}B\geqslant 4$ and $B\not\subset\mathop{\rm
Sing}F$,

(3) a point of general position $o\in B$ is a quadratic
singularity of rank $\geqslant 7$ on $F$ (for $M\geqslant 7$),

(4) a general point $o\in B$ is a quadratic singularity of rank
$\in\{3,4,5,6\}$ on the hypersurface $F$.

We have to exclude each of these four options.

{\bf Proposition 2.1.} {\it The option (1) does not take place:
the inequality} $\mathop{\rm codim}B\geqslant 4$ {\it holds.}

{\bf Proof.} Assume that (1) takes place. Then for $M\geqslant 6$
and for $M=5$ when $\mathop{\rm codim}B=2$ take any irreducible
curve $C\subset B$, such that $C\cap\mathop{\rm Sing}F=\emptyset$.
It is well known (see, for instance, \cite[Chapter 2, Lemma 2.1
and Proposition 2.3]{Pukh13a}), that for every divisor $D\sim nH$
the inequality $\mathop{\rm mult}_CD\leqslant n$ holds. Since $B$
is the centre of a maximal singularity of the system $\Sigma$ and
$B\not\subset\mathop{\rm Sing}F$, we get $\mathop{\rm mult}_BD>n$,
so that $\mathop{\rm mult}_CD>n$, either, for a general divisor
$D\in\Sigma$. We obtained a contradiction that proves our claim
for $M\geqslant 6$ and for $M=5$ in the case $\mathop{\rm
codim}B=2$.

Let us consider the remaining case $M=5$ and $\mathop{\rm
codim}B=3$, that is, $B\subset F$ is an irreducible curve on the
four-dimensional quintic $F\subset{\mathbb P}^5$. If $B$ does not
contain singular points of the hypersurface $F$, one can argue as
above, however, if $B\cap\mathop{\rm Sing}F\neq\emptyset$, then
these arguments do not work. Let $Z=(D_1\circ D_2)$ be the
self-intersection of the system $\Sigma$, that is, the algebraic
cycle of the scheme-theoretic intersection of general divisors
$D_1,D_2\in\Sigma$. Since $B\subset F$ is a curve, the inequality
$\mathop{\rm mult}_BZ>4n^2$ holds. Furthermore, the degree of the
2-cycle $Z$ in ${\mathbb P}^5$ is $5n^2$. Therefore, there is a
surface $S$ on $F$, satisfying the inequality
\begin{equation}\label{23.11.23.1}
\mathop{\rm mult}\nolimits_BS>\frac45\mathop{\rm deg}S.
\end{equation}

Consider first the case when $B$ is a line in ${\mathbb P}^5$. Let
$p\in B$ be an arbitrary point, non-singular on $F$. If
$$
S\not\subset T_pF,
$$
then the intersection $(S\circ T_pF)$ is an effective 1-cycle in
${\mathbb P}^5$ of degree $5n^2$, satisfying the inequality
$\mathop{\rm mult}_p(S\circ T_pF)>8n^2$, which is impossible.
Therefore,
$$
S\subset\bigcap_{p\in B}T_pF.
$$

However, for points $s,p\in B$ of general position the tangent
hyperplanes $T_sF$ and $T_pF$ are distinct, so that $T_sF\cap
T_pF\cap F$ is the section of $F$ by some linear subspace of
dimension 3 and by assumption a general point $q\in B$ is
non-singular on that section. Therefore, $S$ is an irreducible
component of the 2-cycle $(T_sF\circ T_pF\circ F)$ of degree 5 in
${\mathbb P}^5$, non-singular at a general point of the line $B$,
but satisfying the inequality (\ref{23.11.23.1}). The only
possibility is that the irreducible surface $S$ is a plane. But
$F$ does not contain planes. Therefore, the curve $B$ can not be a
line.

Let $p,q\in B$ be distinct points, non-singular on $F$. Since
$$
\mathop{\rm mult}\nolimits_pS+\mathop{\rm
mult}\nolimits_qS>\frac85\mathop{\rm deg}S,
$$
the line $[p,q]\subset{\mathbb P}^5$, joining $p$ and $q$, is
contained in $S$, and therefore in $F$. If $B$ is a curve in some
2-plane, then $S$ is that plane $\langle B\rangle$, which is
contained in $F$, which contradicts the assumption. Therefore, the
linear span $\langle B\rangle$ is not a plane, but then the secant
set
$$
\mathop{\rm Sec}(B)=\overline{\bigcup_{p,q\in B}[p,q]}
$$
is three-dimensional, which is again impossible, as $\mathop{\rm
Sec}(B)\subset S$. (The union is taken over all pairs of distinct
points on $B$.) Proof of Proposition 2.1 is complete.

Therefore, one of the following options takes place: (2), (3) or
(4). Let $o\in B$ be a point of general position.

{\bf Proposition 2.2.} {\it Assume that $o\not\in\mathop{\rm
Sing}F$ (that is, the option (2) takes place). Then there is an
irreducible subvariety $Y\subset F$ of codimension 3, satisfying
the inequality}
\begin{equation}\label{24.11.23.1}
\mathop{\rm mult}\nolimits_oY>\frac{8}{M}\mathop{\rm deg}Y.
\end{equation}

{\bf Proof:} see Subsection 2.3.

{\bf Proposition 2.3.} {\it Assume that the option (3) takes
place. Then there is an irreducible subvariety $Y\subset F$ of
codimension 2, satisfying the inequality (\ref{24.11.23.1}).}

{\bf Proof:} see Subsection 2.3.

{\bf Proposition 2.4.} {\it Assume that the option (4) takes
place. Then there is an irreducible subvariety $Y\subset F$ of
codimension 2, satisfying the inequality}
$$
\mathop{\rm mult}\nolimits_oY>\frac{4}{M}\mathop{\rm deg}Y.
$$

{\bf Proof:} see Subsection 2.3.

Let show now that on the variety $F$ there can be no subvariety
$Y$, the existence of which is claimed by propositions 2.2, 2.3
and 2.4.\vspace{0.3cm}


{\bf 2.2. Hypertangent divisors.} For a fixed point $o\in F$ and
$j\in\{2,\dots,M-1\}$ by the symbol $\Lambda_j$ we denote the
$j$-th {\it hypertangent linear system} at that point, see
\cite[Chapter 3]{Pukh13a}. In the coordinate notations of
Subsection 0.1, where $o=(0,\dots,0)$, on the intersection of $F$
with the affine chart ${\mathbb A}^M_{z_1,\dots,z_M}$ the linear
system $\Lambda_j$ has the form
$$
\left|\left.\sum^j_{i=1}s_{j,j-i}(z_*)f_{[1,i]}(z_*)\right|_F=0\right|,
$$
where $f_{[1,i]}=q_1+\dots+q_i$ is the left segment of the
polynomial $f$ of length $i$ and the homogeneous polynomials
$s_{j,j-i}$ in $z_1,\dots,z_M$ run through the spaces of
homogeneous polynomials ${\cal P}_{j-i,M}$ of degree $j-i$
independently of each other.

The technique of hypertangent divisors and the way it is applied
for proving the birational rigidity are well known, see
\cite[Chapter 3]{Pukh13a}; we will only note the key points of the
arguments in the assumptions of Propositions 2.2-2.4.

{\bf Exclusion of the case (2).} Assume that this case takes
place. Since $\mathop{\rm mult}_o Y\leqslant\mathop{\rm deg}Y$ for
every irreducible subvariety, the inequality (\ref{24.11.23.1})
implies that $M\geqslant 9$. The condition (R1) implies that for
$j\in\{5,6,\dots,M-1\}$ the inequality
$$
\mathop{\rm codim}\nolimits_o(\mathop{\rm Bs}\Lambda_j\subset
F)\geqslant j-4
$$
holds, so that we can construct a sequence of irreducible
subvarieties
$$
Y_3=Y,\, Y_4,\,\dots,\, Y_{M-5}
$$
(in the notations of Proposition 2.2), where $\mathop{\rm
codim}(Y_i\subset F)=i$ and for $i\in\{3,\dots,M-6\}$ the
subvariety $Y_{i+1}$ is an irreducible component of the effective
cycle $(Y_i\circ D_{i+5})$, where $D_{i+5}\in\Lambda_{i+5}$ is a
general divisor, with the maximal value of the ratio $\mathop{\rm
mult}_o\slash\mathop{\rm deg}$ among all components of that cycle.
This construction makes sense, because $\mathop{\rm
codim}_o\mathop{\rm Bs}\Lambda_{i+5}\geqslant i+1$, so that $Y_i$
is not contained in the support of the divisor $D_{i+5}$. For the
last subvariety $Y_{M-5}$ in this sequence we get the inequality
$$
\mathop{\rm mult}\nolimits_oY_{M-5}>\left(\frac{8}{M}\cdot\frac98
\cdot\dots\cdot\frac{M}{M-1}\right) \mathop{\rm deg}Y_{M-5}.
$$
The expression in brackets is equal to 1, which gives a
contradiction, excluding the case (2).

{\bf Exclusion of the case (3).} Assume that this case takes
place. Here again $M\geqslant 9$. Our arguments are completely
similar to the arguments in the case (2) with the obvious
modifications, which we will point out. Here for
$j\in\{6,\dots,M-1\}$ by (R2) the inequality
$$
\mathop{\rm codim}\nolimits_o(\mathop{\rm Bs}\Lambda_j\subset
F)\geqslant j-5
$$
holds, so that we construct (in exactly the same way) a sequence
of irreducible subvarieties
$$
Y_2=Y,\, Y_3,\,\dots,\, Y_{M-6},
$$
where $\mathop{\rm codim}(Y_i\subset F)=i$ and the last subvariety
$Y_{M-6}$ satisfies the inequality
$$
\mathop{\rm mult}\nolimits_oY_{M-6}>\left(\frac{8}{M}\cdot\frac98
\cdot\dots\cdot\frac{M}{M-1}\right)\mathop{\rm deg}Y_{M-6},
$$
which is impossible: the case (3) is excluded.

{\bf Exclusion of the case (4).} Assume that this case takes
place. Here $M\geqslant 5$ and for $j\in\{2,\dots,M-1\}$ by (R3)
the equality
$$
\mathop{\rm codim}\nolimits_o(\mathop{\rm Bs}\Lambda_j\subset
F)=j-1
$$
holds, so that we construct a sequence of irreducible subvarieties
$$
Y_2=Y,\, Y_3,\,\dots,\, Y_{M-2},
$$
where the last subvariety (it is an irreducible curve) $Y_{M-2}$
satisfies the inequality
$$
\mathop{\rm mult}\nolimits_oY_{M-2}>\left(\frac{4}{M}\cdot\frac54
\cdot\dots\cdot\frac{M}{M-1}\right)\mathop{\rm deg}Y_{M-2},
$$
which gives the required contradiction, excluding the case (4).

Proof of the birational rigidity of the hypersurface $F$ is
completed. Q.E.D.\vspace{0.3cm}


{\bf 2.3. Local inequalities.} Let us consider again the
self-intersection $Z=(D_1\circ D_2)$ of the linear system
$\Sigma$, where $D_1,D_2\in\Sigma$ is a general pair of divisors.
In the assumptions of Proposition 2.2 let $P\subset F$ be the
section of $F$ by a general $\mathop{\rm codim}(B\subset{\mathbb
P}^M)$-dimensional linear subspace in ${\mathbb P}^M$, containing
the point $o\in B$. Since $\mathop{\rm codim}(B\subset{\mathbb
P}^M)\geqslant 5$, we have $\dim P\geqslant 4$ and may assume that
$o\in P$ is an isolated centre of a non-canonical singularity of
the pair $(P,\frac{1}{n}D_P)$, where $D_P\in\Sigma_P$ is a general
divisor and $\Sigma_P=\Sigma|_P$. It is well known (see, for
instance, \cite[Chapter 2, Theorem 4.1]{Pukh13a} or \cite{Ch05c}),
that in this case the $8n^2$-inequality holds: for some linear
subspace $\Theta(P)\subset E_P$ of codimension 2, where
$E_P\subset P^+$ is the exceptional divisor of the blow up $P^+\to
P$ of the point $o$ on $P$, the estimate
$$
\mathop{\rm mult}\nolimits_o Z_P+\mathop{\rm
mult}\nolimits_{\Theta(P)} Z^+_P>8n^2
$$
holds, where $Z_P$ is the self-intersection of the linear system
$\Sigma_P$ and $Z^+_P$ is its strict transform on $P^+$. Moreover,
if $\mathop{\rm mult}\nolimits_o Z_P\leqslant 8n^2$, then the
subspace $\Theta(P)$ is uniquely determined by the pair
$(P,\frac{1}{n}D_P)$ (that is, by the system $\Sigma_P$). Coming
back to the original variety $F$, we see that the
$8n^2$-inequality holds already on $F$: for some linear subspace
$\Theta\subset E_F$ of codimension 2, where $E_F\subset F^+$ is
the exceptional divisor of the blow up $F^+\to F$ of the point $o$
on $F$, the estimate
$$
\mathop{\rm mult}\nolimits_o Z+\mathop{\rm mult}\nolimits_{\Theta}
Z^+>8n^2
$$
holds, where $Z^+$ is the strict transform of $Z$ on $F^+$ and, of
course, $\Theta\cap P^+=\Theta(P)$. If $R\subset F$ is a general
hyperplane section, containing the point $o$, such that
$R^+\supset \Theta$, then it is easy to see that
$$
\mathop{\rm mult}\nolimits_o Z_R>8n^2=\frac{8}{M} \mathop{\rm deg}
Z_R,
$$
where $Z_R=(Z\circ R)=Z|_R$ is the self-intersection of the mobile
linear system $\Sigma|_R$. By linearity of the last inequality
there is an irreducible component $Y$ of the effective cycle
$Z_R$, satisfying the inequality (\ref{24.11.23.1}), which
completes the proof of Proposition 2.2.

Proposition 2.3 follows directly from the generalized
$4n^2$-inequality \cite{Pukh2017b} (see also \cite[Chapter I, \S
2]{Pukh2022b}).

Now let us prove Proposition 2.4. In \cite{Pukh2017a} the
following general fact was shown.

{\bf Proposition 2.5.} {\it Let $X$ be a variety with quadratic
singularities of rank at least 4, and assume that $\mathop{\rm
codim}(\mathop{\rm Sing}X\subset X)\geqslant 4$. Assume further
that some divisor $E$ over $X$ is a non-canonical singularity of
the pair $(X,\frac{1}{n}\Sigma)$ with the centre
$B\subset\mathop{\rm Sing}X$, where $\Sigma$ is a mobile linear
system. Then the self-intersection $Z$ of the system $\Sigma$
satisfies the inequality}
$$
\mathop{\rm mult}\nolimits_BZ>4n^2.
$$
(This is \cite[\S 2, Proposition 2.4]{Pukh2017a}.)

This fact was shown for quadratic singularities of rank $\geqslant
5$ in \cite[Section 3]{EcklPukhlikov2014}, however, there was an
incorrectness in the proof, corrected in \cite[\S 2]{Pukh2017a}:
applying the technique of counting multiplicities, one should use,
instead of the numbers of paths $p_{ij}$ in the graph $\Gamma$ of
the maximal singularity, the coefficients $r_{ij}$ that take into
account the singularities of the exceptional divisors. The
relation between the two groups of integral coefficients is
described in detail in \cite[\S 2]{Pukh2017a}.

If the rank of the quadratic singularity $o\in F$ is at least 4,
then we apply Proposition 2.5 and get:
$$
\mathop{\rm mult}\nolimits_oZ>4n^2=\frac{4}{M}\mathop{\rm deg} Z,
$$
which immediately implies Proposition 2.4. Therefore, it remains
to consider the only case when $B=\{o\}$ is a quadratic
singularity of rank 3. Let
$$
\varphi_{j,j-1}\colon F_j\to F_{j-1}
$$
be the resolution of the maximal singularity $E$, that is, the
sequence of blow ups of irreducible subvarieties $B_{j-1}\subset
F_{j-1}$, where $j=1, \dots, N$, $F_0=F$, $B_0=\{o\}$ and for
$j=2, \dots, N$ the subvariety $B_{j-1}$ of codimension $\geqslant
2$ is the centre of $E$ on $F_{j-1}$, and the exceptional divisor
$E_N\subset F_N$ of the last blow up is the centre of $E$ on
$F_N$. Set $E_j=\varphi^{-1}_{j,j-1} (B_{j-1})$ for $j=1, \dots,
N$. In order for the proof of the $4n^2$-inequality, given in
\cite[Section 3]{EcklPukhlikov2014} (in the present case, it has
the form $\mathop{\rm mult}\nolimits_oZ>4n^2$) to work for maximal
singularity, the centre of which is the point $o$ of rank 3, the
following properties need to be satisfied:

--- the variety $F_{j-1}$ should be factorial at a general point
of the subvariety $B_{j-1}$, $j=1, \dots, N$,

--- if $\mathop{\rm codim} (B_{j-1}\subset F_{j-1})\geqslant 3$,
then for the discrepancy we should have $a(E_j,F_{j-1})\geqslant
2$, where $j\in\{1, \dots, N\}$.

For $j=1$ these conditions are satisfied by our assumptions about
the hypersurface $F$ (obviously, $a(E_1,F)=M-3\geqslant 2$).
Furthermore, as it was shown in Subsection 1.3, the variety $F_1$
in a neighborhood of the exceptional divisor $E_1$ has at most
quadratic singularities of rank $\geqslant 4$ and $\mathop{\rm
codim} (\mathop{\rm Sing} F_1\subset F_1)\geqslant 4$. Therefore,
the two conditions given above are satisfied for $j\in\{2, \dots,
N\}$, too. This completes the proof of Proposition 2.4.


\section{Codimension of the complement}

In this section we prove Theorem 0.3. We consider the violation of
each of the conditions defining the subset ${\cal F}\subset {\cal
P}$, and estimate the codimension, corresponding to that
violation.\vspace{0.3cm}

{\bf 3.1. The complement ${\cal P}\setminus {\cal F}$.} The set
${\cal P}\setminus {\cal F}$ consists of polynomials $f$, such
that the hypersurface $F(f)$ does not satisfy one of the
conditions of general position, listed in Subsection 0.1. In order
to prove Theorem 0.3, we have to check that violation of each of
these conditions imposes on $f$ at least $\gamma(M)$ independent
conditions.

First of all, \cite[Theorem 3.1]{Pukh2017a} gives us that the set
of polynomials $f\in{\cal P}$, such that $\mathop{\rm Sing} F(f)$
is positive-dimensional, is of codimension at least
$M(M-2)>\gamma(M)$ in ${\cal P}$. Therefore we may assume that the
set $\mathop{\rm Sing} F(f)$ is finite for all polynomials $f$
under consideration. (We did not include the requirement that the
$\mathop{\rm Sing} F(f)$ is zero-dimensional into the list of
conditions of general position in Subsection 0.1 in order to state
Theorems 0.1 and 0.2 in the maximal generality.)

Furthermore, it is easy to check that ${M-1\choose 2}+1$ is
precisely the codimension of the set of polynomials $f\in{\cal
P}$, such that for some point $o\in F(f)$ this point is either a
quadratic singularity of rank $\leqslant 2$ of the hypersurface
$F(f)$, or a singular point of multiplicity $\geqslant 3$.
Therefore, we can consider only the polynomials $f\in{\cal P}$,
such that the set $\mathop{\rm Sing} F(f)$ is finite and consists
of quadratic singularities of rank $\geqslant 3$. Denote the set
of such polynomials by the symbol ${\cal P}^*$. As have seen
above, $\mathop{\rm codim}(({\cal P}\setminus {\cal P}^*)\subset
{\cal P})\geqslant \gamma(M)$.

For a point $o\in{\mathbb P}^M$ let ${\cal P}^*(o)\subset {\cal
P}^*$ be the subset (of codimension 1), consisting of $f\in {\cal
P}^*$, such that $f(o)=0$.

Denote by the symbols ${\cal B}_G$, ${\cal B}_1$, ${\cal B}_2$,
${\cal B}_3$ the subsets in ${\cal P}^*(o)$, consisting of $f\in
{\cal P}^*(o)$, such that, respectively:

--- the point $o$ is a quadratic singularity of rank 3 and the
condition (G) is not satisfied,

--- the point $o$ is non-singular on $F(f)$ and the condition
(R1) is not satisfied,

--- the point $o\in F(f)$ is a quadratic singularity of rank
$\geqslant 7$ (where $M\geqslant 7$) and the condition (R2) is not
satisfied,

--- the point $o\in F(f)$ is a quadratic singularity of rank
$\in \{3,4,5,6\}$ and the condition (R3) is not satisfied.

In order to prove Theorem 0.3, it is sufficient to check (taking
into account that the point $o$ varies in ${\mathbb P}^M$) that
the codimension of each of these four subsets in ${\cal P}^*(o)$
is at least $\gamma(M)+M-1$. (It is easy to see that the set of
four-dimensional quintics that either contain a two-dimensional
plane or have a section by a three-dimensional subspace in
${\mathbb P}^5$ with a whole line of singular points, is of
codimension higher than $\gamma(5)=6$.)\vspace{0.3cm}


{\bf 3.2. Quadratic points of rank 3.} Let us estimate he
codimension of the set ${\cal B}_G$ in ${\cal P}^*(o)$. For any
$M\geqslant 5$ the condition that the point $o$ is a singularity
of the hypersurface $F(f)$ means in the notations of Subsection
0.1 that $q_1\equiv 0$ ($M$ independent conditions), and the
condition that the rank of the quadratic form $q_2$ is 3 gives
${M-2\choose 2}$ additional independent conditions. We apply again
the estimate of \cite[Theorem 3.1]{Pukh2017a}, this time to the
cubic hypersurface $\{q_3(0,0,0,z_4,\dots,z_M)=0\}$ in ${\mathbb
P}^{M-4}$ for $M\geqslant 8$: we get $3(M-6)$ additional
independent conditions, if the singular locus of this hypersurface
is positive-dimensional. Finally, if that singular set is
zero-dimensional, then the violation of the last part of the
condition (G) gives an additional codimension 1. Therefore, in the
case of the set ${\cal B}_G$ for $M\geqslant 8$ it is sufficient
to check the inequality
$$
{M-2\choose 2}+3M-17\geqslant {M-1\choose 2},
$$
which is very easy. For $M\in\{5,6,7\}$ the inequality
$$
\mathop{\rm codim}({\cal B}_G\subset {\cal P}^*(o))\geqslant
\gamma(M)+M-1
$$
is easy to check, too.\vspace{0.3cm}


{\bf 3.3. The regularity condition at a non-singular point.}
Assume that $M\geqslant 6$. The set ${\cal B}_1$ consists of
polynomials $f$, such that the sequence (\ref{05.12.23.1}) is not
regular. We apply the well known method of estimating the
codimension of that set (see \cite{Pukh98b} or \cite[Chapter
3]{Pukh13a} or \cite{EvansPukh2019}): for each $a\in\{6,\dots,M\}$
consider the subset ${\cal B}_{1,a}$, consisting of polynomials
$f$, such that the regularity of the sequence (\ref{05.12.23.1})
is violated {\it for the first time} by the polynomial
$q_a|_{T_oF}$, that is, the set of common zeros of the system of
polynomials $q_i$, where $i\in\{6,\dots,a-1\}$, has the correct
codimension, whereas $q_a$ vanishes on one of the components of
that set. Since the polynomials $q_i$ are homogeneous, it is
natural to consider the projectivization ${\mathbb P}(T_oF)\cong
{\mathbb P}^{M-2}$. Now we have
$$
{\cal B}_1={\cal B}_{1,6}\bigsqcup\cdots\bigsqcup {\cal B}_{1,M}
$$
and the ``projection method'' (see any of the references above)
estimates the codimension of the subsets ${\cal B}_{1,a}$ from
below by the integers, respectively,
$$
{M+4\choose 6},\, \dots,\,{M+4\choose a},\,\dots,\, {M+4\choose
M}={M+4\choose 4}.
$$
It is easy to see that the minimum of that set of integers is
${M+4\choose 4}$. Therefore, the codimension of ${\cal B}_1$ is
not less than that number, which is much higher than what we need.
\vspace{0.3cm}


{\bf 3.4. Quadratic points of higher rank.} Here $\mathop{\rm rk}
q_2\geqslant 7$ by assumption, so that $q_2\not\equiv 0$, however,
$q_1\equiv 0$, which gives $M$ independent conditions. We estimate
the codimension of the set ${\cal B}_2$, consisting of polynomials
$f$, such that the sequence (\ref{05.12.23.2}) is not regular, in
the same way as in Subsection 3.3: ${\cal B}_2={\cal
B}_{2,7}\sqcup\cdots\sqcup{\cal B}_{2,M}$ and the codimension of
the subset ${\cal B}_{2,a}$ is at least
$$
{M+5\choose a},
$$
where $a\in\{7,\dots,M\}$. The minimum of these integers is
${M+5\choose 5}$: this gives an estimate from below for the
codimension of the set ${\cal B}_2$, which is much higher than
what we need:
$$
\mathop{\rm codim} ({\cal B}_2\subset {\cal P}^*(o))\geqslant
{M+5\choose 5}+M.
$$


{\bf 3.5. Quadratic points of lower rank.} Here $\mathop{\rm rk}
q_2\leqslant 6$, which gives the codimension ${M-5\choose 2}$ (for
$M\in\{5,6\}$ the value of this expression is set to be zero),
besides, $q_1\equiv 0$, which gives the additional codimension
$M$. At the same time, $\mathop{\rm rk} q_2\geqslant 3$, so that
$q_2\not\equiv 0$ and the sequence (\ref{05.12.23.3}) is regular
in the first term. As in Subsections 3.3 and 3.4, we present the
set ${\cal B}_{3}$ as a union of disjoint subsets
$$
{\cal B}_{3}={\cal B}_{3,3}\bigsqcup {\cal B}_{3,4}\bigsqcup\cdots
\bigsqcup{\cal B}_{3,M},
$$
where ${\cal B}_{3,a}$ consists of polynomials $f$, such that the
regularity of the sequence (\ref{05.12.23.3}) is violated {\it for
the first time} in $q_a$. Using the same ``projection method'', we
get the estimate
$$
\mathop{\rm codim} {\cal B}_{3,a}\geqslant {M-5\choose 2} +
{M+1\choose a} + M,
$$
$a\in\{3,4,\dots,M\}$. For $a\leqslant M-1$ the estimate is strong
enough for our purposes. However, for $a=M$ this is not the case,
and we have to estimate the codimension of the set ${\cal
B}_{3,M}$ separately. We consider $q_i$, $2\leqslant i\leqslant
M$, as homogeneous polynomials on ${\mathbb P}^{M-1}$. By the
definition of the set ${\cal B}_{3,M}$, if $f\in {\cal B}_{3,M}$,
then the set of common zeros of the polynomials $q_i$, $2\leqslant
i\leqslant M-1$, in ${\mathbb P}^{M-1}$ is zero-dimensional and
for some irreducible curve $C\subset{\mathbb P}^{M-1}$, which is a
component of that set, we have $q_M|_C\equiv 0$. Setting
$\dim\langle C\rangle=b\in\{1,\dots,M-1\}$, we get the
presentation
$$
{\cal B}_{3,M}={\cal B}_{3,M,1}\bigcup{\cal B}_{3,M,2}\bigcup
\dots \bigcup{\cal B}_{3,M,M-1},
$$
where ${\cal B}_{3,M,b}$ consists of $f\in{\cal B}_{3,M}$, such
that $\dim\langle C\rangle=b$. (This union, generally speaking, is
not disjoint.) It is sufficient to estimate from below the
codimension of each set ${\cal B}_{3,M,b}$.

Elementary computations show that
$$
\mathop{\rm codim} {\cal B}_{3,M,1}\geqslant {M-5\choose
2}+{M\choose 2}+M+2.
$$
This is much better than we need. If $b=2$, that is,
$C\subset\langle C\rangle\cong{\mathbb P}^2$ is a plane curve of
degree $d_C\geqslant 2$, then elementary computations show that
the condition $q_i|_C\equiv 0$ for all $i\in\{2,\dots,M\}$ imposes
on the tuple $\{q_i\}$, that is, on the polynomial $f$,
$$
M^2-M+1
$$
conditions for $d_C=2$ (and arbitrary plane $\langle C\rangle$),
and
$$
\sum^{d_C-1}_{i=2}{i+2\choose 2}+ \sum^M_{i=d_C+1}\left(
{i+2\choose 2}-{i+2-d_C\choose 2}\right)
$$
independent conditions for $d_C\in\{3,\dots,M-1\}$ for a fixed
plane $\langle C\rangle\subset{\mathbb P}^{M-1}$, which, taking
into account the variation of that plane gives a stronger estimate
for the codimension than for $d_C=2$, so that in the end we get
the inequality
$$
\mathop{\rm codim} {\cal B}_{3,M,2}\geqslant {M-5\choose 2}+M^2+1,
$$
which is again much better than what we need.

Assume that $b\geqslant 3$. Fix a subspace $P\subset{\mathbb
P}^{M-1}$ of dimension $b$ and consider the subset ${\cal
B}_{3,M,b}(P)$, given by the condition $\langle C\rangle=P$. Here
we will need the technique of {\it good  sequences and associated
subvarieties}, introduced in \cite{Pukh01} (see also \cite[Chapter
3, Section 3]{Pukh13a} or \cite[Definition 3.1]{Pukh2017a}). Let
us remind briefly this approach. Although the curve $C$ is an
irreducible component of the set of common zeros of the
polynomials $q_i$, $i\in\{2,\dots,M-1\}$, for $b\leqslant M-2$,
generally speaking, it is not true that one can choose from these
polynomials $(b-1)$ ones, such that $C$ is a component of their
restrictions onto $P$. However, it is true that one can choose
$(b-1)$ polynomials $q_{i_1},\dots,q_{i_{b-1}}$, where
$i_{\alpha}\in\{2,\dots,M-1\}$, such that there is a sequence of
irreducible subvarieties $\Delta_0=P$,
$\Delta_1,\dots,\Delta_{b-1}$ of the subspace $P$, such that
$\mathop{\rm codim}(\Delta_{\alpha}\subset P)=\alpha$,
$\Delta_{\alpha+1}\subset\Delta_{\alpha}$,
$q_{i_{\alpha+1}}|_{\Delta_{\alpha}}\not\equiv 0$ for
$\alpha=0,\dots,b-2$ and $\Delta_{\alpha+1}$ is an irreducible
component of the closed set
$\{q_{i_{\alpha+1}}|_{\Delta_{\alpha}}=0\}$, and moreover,
$\Delta_{b-1}=C$. For a fixed tuple $q_{i_1},\dots,q_{i_{b-1}}$
consider the set of all sequence $\Delta_0,\dots,\Delta_{b-1}$,
described above, with a free end (that is, we remove the last
condition $\Delta_{b-1}=C$). However, this set is finite, so that
the set of curves, which are the last subvariety $\Delta_{b-1}$ in
one of such sequences, is finite, too (and $C$ is one of these
curves). Thus, fixing the polynomials $q_{i_1},\dots,q_{i_{b-1}}$,
we may assume that the curve $C$ is fixed, either. For
$$
j\not\in\{i_1,\dots,i_{b-1}\}
$$
the condition $q_j|_C\equiv 0$ imposes on $q_j$ at least $jb+1$
independent conditions (since $\langle C\rangle=P$, no polynomial
of degree $j$, which is a product of $j$ linear forms on $P$, can
vanish identically on $C$, see the details in the references
above). The worst estimate for the codimension is obtained when
$\{i_1,\dots,i_{b-1}\}=\{M-b+1,\dots,M-1\}$. Therefore, for
$b=M-1$ we have $M(M-1)+1$ independent conditions (since
$q_M|_C\equiv 0$), and for $3\leqslant b\leqslant M-2$ we have
$$
b(M+2+\dots+M-b)+(M-b)=b(M+\frac12(M-b+2)(M-b-1))+(M-b)
$$
independent conditions. Taking into account the variation of the
subspace $P$ in ${\mathbb P}^{M-1}$, we obtain the estimate for
the codimension as the integer $h(b)$, where
$b\in\{3,\dots,M-1\}$, and
$$
h(t)=\frac12[t^3+(-2M+1)t^2+(M^2+M)t+2].
$$
Considering the derivative
$$
h'(t)=\frac12[3t^2+2(-2M+1)t+(M^2+M)],
$$
it is easy to check that for $M\in\{5,6\}$ the polynomial $h(t)$
is increasing on the segment $[3,M-1]$, so that the worst estimate
for the codimension is given by the number $h(3)$. It is easy to
see that this estimate is strong enough for the proof of Theorem
0.3 for $M=5,6$. If $M\geqslant 7$, then the behaviour of the
function $h(t)$ on the segment $[3,M-1]$ is more complicated: the
both roots $t_*<t^*$ of the derivative $h'(t)$ lie on that
segment, so that $h(t)$ increases on the segment $[3,t_*]$,
decreases on the segment $[t_*,t^*]$ and then increases on
$[t^*,\infty)$. It is easy to check that $M-2\leqslant
t^*\leqslant M-1$, so that the worse estimate from below for the
codimension is the minimum of the following three integers:
$$
h(3)=\frac32 M(M-5)+19,\quad h(M-2)=M(M-1)-1
$$
and $h(M-1)=M(M-1)+1$. For $M=7$ the minimum is $h(M-2)$, for
$M\geqslant 8$ the minimum is $h(3)$. It is easy to see that this
estimate is strong enough. For instance, for $M\geqslant 8$ the
codimension is at least
$$
{M-5\choose 2}+\frac32M(M-5)+M+19,
$$
which is much higher than what we need.

Proof of Theorem 0.3 is complete. Q.E.D.

Note that for $M\geqslant 8$ the estimate for the codimension of
the complement ${\cal P}\setminus {\cal F}$ in Theorem 0.3 is
related only to the violation of the condition for the rank of the
quadratic singularity; violation of the condition (G) or the
regularity conditions (R1-3) gives a higher codimension.


\begin{flushleft}
Department of Mathematical Sciences,\\
The University of Liverpool
\end{flushleft}

\noindent{\it pukh@liverpool.ac.uk}

\end{document}